\documentclass[12pt,twoside]{article}
\usepackage{amsmath,amsbsy,amsfonts}
\usepackage{graphicx,color}
\newtheorem{theorem}{Theorem}[section]
\newtheorem{lemma}[theorem]{Lemma}
\newtheorem{proposition}[theorem]{Proposition}
\newtheorem{corollary}[theorem]{Corollary}
\newtheorem{remark}[theorem]{Remark}
\newtheorem{claim}[theorem]{Claim}

\newcommand{\fim}{\hfill\rule{2mm}{2mm}}

\newcommand{\ds}{\displaystyle}
\newcommand{\R}{\mathbb{R}}

\begin{document}

\setlength{\baselineskip}{4.5mm} \setlength{\oddsidemargin}{8mm}
\setlength{\topmargin}{-3mm}
\title{\Large\sf Existence and nonexistence of least energy nodal solution for a class of elliptic problem in $\mathbb{R}^{2}$ }
\author{\sf
 Claudianor O. Alves\thanks{Partially supported by CNPq - Grant 304036/2013-7, coalves@dme.ufcg.edu.br} \;\; and \;\;  Denilson S. Pereira\thanks{denilsonsp@dme.ufcg.edu.br} \\
 Universidade Federal de Campina Grande\\
 Unidade Acad\^emica de Matem\'atica - UAMat\\
 CEP: 58.429-900 - Campina Grande - PB - Brazil\\ }

\pretolerance10000
\date{}
\numberwithin{equation}{section} \maketitle
\begin{abstract}
In this work, we prove the existence of least energy nodal solution for a class of elliptic problem in both cases, bounded and
unbounded domain, when the nonlinearity has exponential critical growth in $\mathbb{R}^2$. Moreover, we also prove a nonexistence result of least energy  nodal solution for the autonomous case in whole $\mathbb{R}^{2}$. \\


\vspace{0.3cm}


\end{abstract}

\section{Introduction}

This paper concerns with the existence of least energy nodal
solution for the following class of elliptic problem
$$
\left\{
\begin{array}{l} -\Delta u+V(x)u=f(u),  \mbox{ in } \Omega, \\
\mbox{}\\
u\in H_0^1(\Omega),
\end{array}
\right.
\eqno (P)
$$
where $\Omega\subset \R^2$ is a smooth bounded domain or $\Omega = \mathbb{R}^2$, $V:\overline{\Omega} \rightarrow\R$
is a continuous function verifying some hypotheses which will be fix later on. Related to the nonlinearity, we assume that $f:\R\rightarrow\R$ is a
$C^1$-function, which can have an exponential critical growth at both $+\infty$ and
$-\infty$, that is, it behaves like $e^{\alpha_0s^2}$, as
$|s|\rightarrow\infty$, for some $\alpha_0>0$. More precisely,
\begin{equation} \label{expg}
\lim_{|s|\rightarrow\infty}\dfrac{f(s)}{e^{\alpha|s|^2}}=0\ \
\forall \alpha>\alpha_0,\ \
\lim_{|s|\rightarrow\infty}\dfrac{f(s)}{e^{\alpha|s|^2}}=\infty\ \
\forall\alpha<\alpha_0 \,\,\, ( \mbox{see} \, \cite{DMR} \, ).
\end{equation}

In the last years, we have observed that the existence of nodal solution has received a special attention of a lot of researches. In Cerami, Solimini and Struwe \cite{CSS}, the authors showed the existence of multiples nodal solutions for the following class of elliptic problem with critical growth
$$
\left\{
\begin{array}{l}
- \Delta{u}-\lambda u=|u|^{2^{*}-2}u, ~~ \mbox{in} ~~ \Omega\\
\mbox{}\\
u=0, ~~ \mbox{on} ~~ \partial \Omega
\end{array}
\right.
\eqno{(P_1)}
$$
where $\Omega=B_R(0) \subset \mathbb{R}^{N}, N \geq 7, 2^*=\frac{2N}{N-2}$ and $\lambda \in [0,\lambda_1]$, with $\lambda_1$ being the first eigenvalue of $(-\Delta, H^{1}_{0}(\Omega))$.  In Bartsch and Willem \cite{BWillem}, infinitely many radial nodal solutions was proved for the problem
$$
\left\{
\begin{array}{l}
- \Delta{u}+u=f(|x|,u), ~~ \mbox{in} ~~ \mathbb{R}^{N}\\
\mbox{}\\
u \in H^{1}(\mathbb{R}^{N}),
\end{array}
\right.
\eqno{(P_2)}
$$
where $f$ is a continuous function with sucritical growth and verifying some hypotheses. In Cao and Noussair \cite{CN}, the authors studied the existence and multiplicity of positive and nodal solutions for the following class of problems
$$
\left\{
\begin{array}{l}
- \Delta{u}+u=Q(x)|u|^{p-2}u, ~~ \mbox{in} ~~ \mathbb{R}^{N}\\
\mbox{}\\
u \in H^{1}(\mathbb{R}^{N})
\end{array}
\right.
\eqno{(P_2)}
$$
by supposing $2<p<\frac{N+2}{N-2}, N \geq 3 $ and some technical conditions on $Q$. In that paper, the main result connects the number of positive and nodal solutions with the number of maximum points of function $Q$.

In Castro, Cossio and Neuberger \cite{CCN} and Bartsch and Wang \cite{BWa1}, the authors studied the existence of nodal solution for a problem like
$$
\left\{
\begin{array}{l}
- \Delta{u}=f(u), ~~ \mbox{in} ~~ \Omega\\
\mbox{}\\
u=0, ~~ \mbox{on} ~~ \partial \Omega
\end{array}
\right.
\eqno{(P_3)}
$$
where $\Omega$ is a smooth bounded domain and $f$ verifies some hypotheses.  In \cite{CCN}, it was assumed that $f$ is superlinear, while that in \cite{BW1}, $f$ is  asymptotically linear at infinity. In Bartsch and Weth \cite{BW1}, existence of multiple nodal solution was also considered for problem $(P_3)$.

In Noussai and Wei \cite{NW1,NW2}, existence and concentration of nodal solution were proved for the problem
$$
\left\{
\begin{array}{l}
- \epsilon^{2}\Delta{u}+u=f(u), ~~ \mbox{in} ~~ \Omega\\
\mbox{}\\
Bu=0, ~~ \mbox{on} ~~ \partial \Omega,
\end{array}
\right.
\eqno{(P_4)}
$$
when $\epsilon \to 0$, where $\Omega$ is smooth bounded domain, $Bu=0$ in \cite{NW1} and $Bu=\frac{\partial u}{\partial \eta}$ in \cite{NW2}.

In Bartsch and Wang \cite{BWa}, the authors have considered the existence and concentration of nodal solution for the following class of problem
$$
\left\{
\begin{array}{l}
- \Delta{u}+(\lambda a(x)+1)u=f(u), ~~ \mbox{in} ~~ \mathbb{R}^{N} \\
\mbox{}\\
u \in H^{1}(\mathbb{R}^{N}),
\end{array}
\right.
\eqno{(P_5)}
$$
when $\lambda \to +\infty$, by supposing that $f$ has a subcritical growth and \linebreak $a:\mathbb{R}^{N} \to \mathbb{R}$ is a nonnegative continuous function with $a^{-1}(\{0\})$ being nonempty and verifying
$$
\mu(\{x \in \mathbb{R}^{N}~;~ a(x) \leq M_0\})<+\infty ~~~~ \mbox{for some} ~~~~~M_0>0.
$$

In \cite{BLW}, Bartsch, Liu and Weth have showed the existence of nodal solution with exactly two nodal regions was established for the problem
$$
\left\{
\begin{array}{l}
- \Delta{u}+a(x)u=f(u), ~~ \mbox{in} ~~ \mathbb{R}^{N}\\
\mbox{}\\
u \in H^{1}(\mathbb{R}^{N}),
\end{array}
\right.
\eqno{(P_6)}
$$
where $a$ is a nonnegative verifying conditions, among which we highlight
$$
\mu(\{x \in B_r(y)~:~ a(x) \leq M\}) \to 0 ~~ \mbox{as} ~~ |y| \to +\infty ~~ \mbox{for any} ~~ M,r>0.
$$

The reader can found more results involving nodal solutions in the papers of Bartsch, Weth and Willem \cite{BWW}, Alves and Soares \cite{ASS}, Bartsch, Clapp and Weth \cite{BCW}, Zou \cite{Z} and their references.

After a literature review, we have observed that there are few papers in the literature where existence of nodal solution has been considered for the case where the nonlinearity has an exponential critical growth. The authors know only the references Adimurthi and Yadava \cite{AY}, Alves \cite{A} and  Alves and Soares \cite{ASS1}. In \cite{AY}, the authors have proved the infinite many radial solution for problem $(P_3)$ when $\Omega=B_R(0) \subset \mathbb{R}^{2}$. In \cite{A}, the authors has proved the existence of nodal solution for a class of problem in exterior domain with Neumann boundary conditions, and in \cite{ASS1}, the existence of nodal solution has been established for a problem like
$$
\left\{
\begin{array}{l}
-\epsilon^{2} \Delta u+ V(x)u=f(u), \,\, \mbox{in} \,\, \mathbb{R}^{N} \\
u \in H^{1}(\mathbb{R}^{N})
\end{array}
\right.
$$
for $\epsilon$ small enough and $V$ verifying some technical conditions.

Motivated by this fact, our goal in the present paper is proving the existence of least energy nodal solution for problem $(P)$ when $\Omega$ is a smooth bounded domain or $\Omega=\mathbb{R}^{2}$.
Here, we also show a nonexistence result of least energy solution for $(P)$ when the potential $V$ is constant.  Once that we will work with exponential critical growth in whole $\mathbb{R}^{2}$,
a key inequality in our arguments is the Trudinger-Moser inequality for bounded domain, see \cite{M} and \cite{T}, which claims that for any  $ u \in H^{1}_0(\Omega)$,
\begin{equation} \label{X0}
\int_{\Omega}
e^{\alpha\left|u\right|^{2}}dx
< +\infty, \,\,\,\, \mbox{ for every }\,\,\alpha >0.
\end{equation}
Moreover, there exists a positive constant $C=C(\alpha,|\Omega|)$ such that
\begin{equation} \label{X1}
\sup_{||u||_{H_0^{1}(\Omega)} \leq 1} \int_{\Omega} e^{\alpha|u|^{2}} dx \leq C , \,\,\,\,\,\,\, \forall \, \alpha  \leq 4 \pi .
\end{equation}
A version of the above inequality in whole space $\R^2$ has been proved by Cao \cite{Cao} and has the following statement:
\begin{equation} \label{X2}
\int_{\R^2}(e^{\alpha|u|^2-1})dx<+\infty,\ \ \mbox{for all} \,\,\, u\in H^1(\R^2)\ \ \mbox{and} \,\, \alpha>0.
\end{equation}
Furthermore, if $\alpha \leq 4\pi$ and $|u|_{L^2(\R^2)}\leq M$, there exists
a constant \linebreak $C_1=C_1(M,\alpha)$ such that
\begin{equation} \label{X4}
\sup_{ |\nabla u|_{L^2(\R^2)}\leq 1} 
\int_{\R^2}(e^{\alpha|u|^2-1})dx\leq C_1.
\end{equation}

Hereafter, the function $f$ satisfies the ensuing assumptions:

\begin{enumerate}

\item[$(f_1)$] There is $C>0$ such that
$$
|f(s)|\leq Ce^{4\pi |s|^2}\ \ \mbox{for all}\ \ s\in\R;
$$

\item[$(f_2)$] $\ds\lim_{s\rightarrow
0}\dfrac{f(s)}{s}=0$;

\item[$(f_3)$] There is $\theta>2$ such that
$$0<\theta F(s):=\theta\int_0^{s}f(t)dt\leq sf(s),\ \ \mbox{for all}\ \
s\in\R\setminus\{0\}.$$

\item[$(f_4)$] The function $s\rightarrow\dfrac{f(s)}{|s|}$ is
strictly increasing in $(0,+\infty)$.

\item[$(f_5)$] There exist constants $p>2$ and $C_p>0$ such  that
$$
sign(s)f(s)\geq C_p|s|^{p-1}\ \ \mbox{for all}\ \ s\in\R\setminus\{0\},
$$
where
$$
sign(s)=\left\{
\begin{array}{l}
\,\,\,\,1, \,\, s >0 \\
-1, \,\, s<0.
\end{array}
\right.
$$
\end{enumerate}

Our main  result related to the case where $\Omega$ is a bounded domain is the following:

\noindent  \begin{theorem}\label{bd}
Let $\Omega$ be a bounded domain and $V:\overline{\Omega} \to \mathbb{R}$ be a nonnegative continuous function. If $(f_1)-(f_5)$ occur, then problem $(P)$
possesses a least energy nodal solution, provided that
$$
C_p>\left[\beta_p\dfrac{2\theta}{\theta-2}\right]^{(p-2)/2},
$$
where
$$
\beta_p=\inf_{\mathcal{M}_{\Omega}}I_p,
$$
$$
\mathcal{M}_{\Omega}=\{u\in H_0^1(\Omega):\ \ u^{\pm}\neq 0 \ \ \mbox{e}\ \ I'_p(u^{\pm})u^{\pm}=0\}
$$
and
$$
I_p(u)=\dfrac{1}{2}\int_{\Omega}\left(|\nabla u|^2+V(x)|u|^2\right)dx-\dfrac{1}{p}\int_{\Omega}|u|^pdx .
$$
\end{theorem}

For the case where $\Omega = \mathbb{R}^{2}$, we have two results. The first one is a nonexistence result of least energy nodal solution whose statement is the following:

\begin{theorem}  \label{ud2}
Suppose that $V(x)=V_0>0$ for all $x\in \mathbb{R}^{2}$ and $f$ satisfies
$(f_1)-(f_5)$. Then, the autonomous problem
$$
\left\{\begin{array}{lcl} -\Delta u+V_0u=f(u),  &\mbox{ in }& \R^2, \\
u\in H^1(\R^2),
\end{array}\right. \leqno (P)
$$
does not have  a least energy nodal solution, provided that
\begin{equation} \label{EQ1}
C_p>\left[\chi _p\dfrac{2\theta}{\theta-2}\right]^{(p-2)/2},
\end{equation}
where
$$
\chi_p=\inf_{\mathcal{M}_{B_1(0)}}J_p,
$$
$$
\mathcal{M}_{B_1(0)}=\{u\in H_0^1(B_1(0)):\ \ u^{\pm}\neq 0 \ \ \mbox{and}\ \ J'_p(u^{\pm})u^{\pm}=0\}
$$
and
$$
J_p(u)=\dfrac{1}{2}\int_{B_1(0)}\left(|\nabla u|^2+V_0 |u|^2\right)dx-\dfrac{1}{p}\int_{B_1(0)}|u|^pdx .
$$

\end{theorem}

Our second result is related to the existence of least energy nodal solution for a non-autonomous problem. For this case, we will assume the ensuing hypotheses on function $V$:

\begin{enumerate}
\item[$(V_1)$] There exists a constant $V_0>0$ such that $V_0\leq V(x)$ for all $x\in\R^2$;

\item[$(V_2)$] There exists a continuous $\mathbb{Z}^2$-periodic function $V_\infty:\R^2\rightarrow\R$
satisfying 
$$
V(x)\leq V_\infty(x) \,\,\, \forall x\in \R^2
$$ 
and
$$
\lim_{|x|\rightarrow\infty}|V(x)-V_\infty(x)|=0.
$$
We recall that a function $h:\mathbb{R}^{2} \to \mathbb{R}$ is $\mathbb{Z}^2$-periodic when
$$
h(x)=h(x+y),\ \ \mbox{for all}\ \ x\in\R^2 \,\,\, \mbox{and} \,\,\, y\in \mathbb{Z}^2.
$$
\item[$(V_3)$] There exist $\mu<1/2$ and $C>0$ such that
$$V(x)\leq V_\infty(x)-Ce^{-\mu|x|},\ \ \mbox{for all}\ \ x\in\R^2.$$
\end{enumerate}

Our main result involving the above hypotheses is the following:

\begin{theorem}\label{ud}
Suppose that hypotheses $(V_1)-(V_3)$, $(f_1)-(f_5)$ and (\ref{EQ1}) are fulfilled. Then the elliptic problem
$$
\left\{\begin{array}{lcl} -\Delta u+V(x)u=f(u),  &\mbox{ in }& \R^2, \\
u\in H^1(\R^2),
\end{array}\right. \leqno (P)
$$
possesses a least energy nodal solution.

\end{theorem}

We conclude this section by giving a sketch of the proofs. The basic
idea goes as follows. To prove Theorem \ref{bd} we will use the
Nehari method and the deformation lemma. Our inspiration comes from of \cite{CSS}, however in that paper the authors used a deformation lemma in cones together with the fact that the nonlinearity is odd. Here, we developed a new approach to get a Palais-Smale sequence of nodal function associated to the least energy nodal level, for details see Section 2. In order to prove Theorem \ref{ud}, we invoke Theorem
\ref{bd} to obtain a sequence $(u_n)$ of least energy nodal
solutions to problem $(P)$ when $\Omega=B_n(0)$. Then, we prove that $(u_n)$ is
weakly convergent in $H^1(\R^2)$, and its weak limit is a least energy nodal
solution of the problem $(P)$.
\section{Bounded Domain}

In this section, we consider the existence of least energy nodal solution for problem $(P)$ when $\Omega$ is a smooth bounded domain. Let us denote by $E$ the Sobolev space $H_0^1(\Omega)$ endowed with the norm
$$
\|u\|^2=\int_{\Omega}\left(|\nabla u|^2+V(x)|u|^2\right)dx.
$$
From assumptions $(f_1)$ and $(f_2)$, given $\epsilon>0$, $q\geq 1$
and $\alpha>4$, there exists a positive constant
$C=C(\epsilon,q,\alpha)$ such that

\begin{equation}\label{1.5}
|sf(s)|,\ \ |F(s)| \leq \epsilon\dfrac{s^2}{2}+C|s|^{q}e^{\alpha \pi s^2},\
\ \mbox{for all}\ \ s\in \R.
\end{equation}
Thus, by Trudinger-Moser inequality (\ref{X0}), $F(u)\in L^1(\R^2)$
for all $u\in E$, from where it follows that Euler-Lagrange
functional associated with $(P)$ $I:E \to \mathbb{R}$ given by
$$
I(u)=\dfrac{1}{2}\|u\|^2-\int_{\Omega}F(u)dx
$$
is well defined. Furthermore, using standard arguments, we see that
$I$ is a $C^1$ functional on $E$ with
$$
I'(u)v=\int_{\Omega}\left[\nabla u\nabla v+V(x)uv\right]dx-\int_{\Omega}f(u)vdx,\ \ \mbox{for all}\ \ v\in E.
$$
Consequently, critical points of $I$ are precisely the weak solution
of problem $(P)$. We know that every nontrivial critical point of
$I$ is contained in the Nehari manifold
$$
\mathcal{N}_\Omega=\{u\in E\setminus\{0\}:\; I'(u)u=0\}.
$$
Since we are interested in least energy nodal solution, we define
the nodal Nehari set
$$
\mathcal{M}_\Omega=\{u\in E:\;u^{\pm}\neq 0,\; I'(u^{\pm})u^{\pm}=0\},
$$
and
$$
c^{*}_\Omega=\inf_{u\in\mathcal{M}_\Omega}I(u).
$$

By a least energy nodal solution, we understand as being a function $u \in \mathcal{M}_\Omega$ such that
$$
I(u)=c^{*}_\Omega ~~ \mbox{and} ~~ I'(u)=0.
$$

Next, we state some necessary results to prove Theorem
\ref{bd}. The proofs of some of them are in Section 4.

\begin{lemma}\label{A}
There exists $A>0$ such that
$$c^{*}_\Omega\leq A <\left(\dfrac{1}{2}-\dfrac{1}{\theta}\right).$$
\end{lemma}

\noindent {\bf Proof.} ~ Let $\tilde{u}\in \mathcal{M}_{\Omega}\subset H_0^1(\Omega)$ verifying
\begin{equation} \label{E1}
I_p(\tilde{u})=\beta_p ~~ \mbox{and} ~~ I'_p(\tilde{u})=0.
\end{equation}
The reader can find the proof of the above claim in Bartsch and Weth \cite{BW}. Once $\tilde{u}^{\pm}\neq 0$ there exist $0<s,t$ such that
$s\tilde{u}^{+},t\tilde{u}^{-}\in\mathcal{N}_\Omega$ and $s\tilde{u}^{+}+t\tilde{u}^{-}\in\mathcal{M}_\Omega$. Then,
$$
c^{*}_\Omega\leq I(s\tilde{u}^{+}+t\tilde{u}^{-})=I(s\tilde{u}^{+})+I(t\tilde{u}^{-}),
$$
loading to
$$
c^{*}_\Omega\leq \dfrac{s^2}{2}\int_{\Omega}\left(|\nabla\tilde{u}^{+}|^2+V(x)|\tilde{u}^{+}|^2\right)dx-\int_{\Omega}F(s\tilde{u}^{+})dx
$$
$$
+\dfrac{t^2}{2}\int_{\Omega}\left(|\nabla\tilde{u}^{-}|^2+V(x)|\tilde{u}^{-}|^2\right)dx-\int_{\Omega}F(t\tilde{u}^{-})dx.
$$
By $(f_5)$,
$$
c^{*}_\Omega\leq \left(\dfrac{s^2}{2}-\dfrac{C_ps^p}{p}\right)\int_{\Omega}|\tilde{u}^{+}|^pdx+\left(\dfrac{t^2}{2}-\dfrac{C_pt^p}{p} \right)\int_{\Omega}|\tilde{u}^{-}|^pdx,
$$
and so,
$$
c^{*}_\Omega\leq\max_{r\geq 0}\left\{\dfrac{r^2}{2}-\dfrac{C_pr^p}{p}\right\}\int_{\Omega}|\tilde{u}|^pdx.
$$
A direct computation gives 
$$
\max_{r\geq 0}\left\{\dfrac{r^2}{2}-\dfrac{C_pr^p}{p}\right\}=C_p^{\frac{2}{2-p}}\left(\dfrac{1}{2}-\dfrac{1}{p}\right),
$$
then 
$$
c^{*}_\Omega\leq C_p^{\frac{2}{2-p}}\left(\dfrac{1}{2}-\dfrac{1}{p}\right)\int_{\Omega}|\tilde{u}|^pdx.
$$
Using $(\ref{E1})$ in the above inequality, we get
\begin{equation} \label{A1}
c^{*}_\Omega\leq C_p^{\frac{2}{2-p}}\beta_p:=A.
\end{equation}
From $(f_5)$,
$$
A<\left(\dfrac{1}{2}-\dfrac{1}{\theta}\right),
$$
finishing the proof. \fim

\vspace{0.5 cm}

The next lemma shows two important limits involving the function.

\begin{lemma}\label{l1.8}

Let $(u_n)$ be a sequence in $E$ satisfying

\begin{enumerate}
\item[(1)] $b:=\displaystyle \sup_{n \in \mathbb{N}}\|u_n\|^2<1;$

\item[(2)] $u_n\rightharpoonup u$ in $H_0^1(\Omega)$ and;

\item[(3)] $u_n(x)\rightarrow u(x)$ a.e. in $\Omega$.
\end{enumerate}

Then,

\begin{equation}\label{1.22}
\lim_{n}\int_{\Omega}f(u_n)u_ndx=\int_{\Omega}f(u)udx
\end{equation}
and
\begin{equation}\label{1.23}
\lim_{n}\int_{\Omega}f(u_n)v dx=\int_{\Omega}f(u)v dx,
\end{equation}
for any $v\in E$.
\end{lemma}

\noindent {\bf Proof.} See Section 4.

\fim

\vspace{0.5 cm}

The below result is very know for problem in $\mathbb{R}^{N}$ with $N \geq 3$. Here, we decide to write its proof, because we are working with exponential critical growth.

\begin{lemma}\label{l1.1}

There exists $m_0 >0$  such that
$$0<m_0\leq\|u\|^2,\ \ \forall u\in\mathcal{N}_\Omega.$$

\end{lemma}

\noindent {\bf Proof.} We start by fixing $q>2$ in $(\ref{1.5})$. Suppose  by
contradiction that above inequality is false. Then, there exists a
sequence $(u_n)\subset\mathcal{N}_\Omega$ such that $\|u_n\|^2\rightarrow
0$, as $n\rightarrow \infty$. Since $u_n\in\mathcal{N}_\Omega$,
$$
\|u_n\|^2=\int_{\Omega}f(u_n)u_ndx.
$$
Then, from $(\ref{1.5})$,
$$
\|u_n\|^2\leq\epsilon |u_n|_2^2+C\int_{\Omega}|u_n|^qe^{\alpha|u_n|^2}dx.
$$
By Sobolev imbedding and H\"{o}lder inequality, 
$$
\|u_n\|^2\leq C_1\epsilon\|u_n\|^2+C|u_n|_{2q}^q\left(\int_{\Omega}e^{2\alpha|u_n|^2}dx\right)^{1/2}.
$$
Using again Sobolev imbedding,
$$
(1-C_1\epsilon)\|u_n\|^2\leq C_2\|u_n\|^q\left(\int_{\Omega}e^{2\alpha|u_n|^2}dx\right)^{1/2}.
$$
Choosing $\epsilon>0$ sufficiently small such that
$C_3:=\dfrac{1-C_1\epsilon}{C_2}>0$, we find that 
\begin{equation}\label{1.8}
0<C_3\leq
\|u_n\|^{q-2}\left(\int_{\Omega}e^{2\alpha|u_n|^2}dx\right)^{1/2}.
\end{equation}
Since $\|u_n\|^2\rightarrow 0$, as $n\rightarrow \infty$, there is
$n_0\in\mathbb{N}$ such that
$$2\alpha\|u_n\|^2\leq 4\pi, \ \ \forall n\geq n_0.$$
From Trudinger-Moser inequality (\ref{X1}), it follows that
$$
\int_{\Omega}e^{2\alpha|u_n|^2}dx=\int_{\Omega}e^{2\alpha\|u_n\|^2\left(\frac{|u_n|}{\|u_n\|}\right)^2}dx\leq\int_{\Omega}e^{4\pi\left(\frac{|u_n|}{\|u_n\|}\right)^2}dx\leq C ~~ \forall n \geq n_0.
$$
Thereby, by $(\ref{1.8})$,
$$
0<\left(\dfrac{C_3}{\sqrt{C}}\right)^{1/(q-2)}\leq \|u_n\|, \ \ \forall n\geq n_0,
$$
which contradicts the fact that $\|u_n\|\rightarrow 0$, as $n\rightarrow \infty$. \fim

\begin{corollary}\label{c1.1}
For all $u\in\mathcal{M}_\Omega$,
$$
0<m_0\leq\|u^{\pm}\|^2\leq\|u\|^2.
$$
\end{corollary}

\begin{corollary}\label{l1.5}
There exists $\delta_2>0$ such that
$$
I(u^{\pm})\geq \delta_2 \,\,\, \mbox{and} \,\,\, I(u)\geq 2\delta_2 ~~ \forall u\in\mathcal{M}_\Omega.
$$

\end{corollary}

\noindent {\bf Proof.} Firstly, observe that if $u\in\mathcal{N}_\Omega$, 
$$
I(u)=I(u)-\dfrac{1}{\theta}I'(u)u=\left( \dfrac{1}{2}- \dfrac{1}{\theta}\right)\|u\|^2-\int_{\Omega}\left(f(u)-\dfrac{1}{\theta}f(u)u\right)dx.
$$
Then, from $(f_3)$ and Lemma \ref{l1.1},
$$
I(u)\geq\left( \dfrac{1}{2}- \dfrac{1}{\theta}\right)\|u\|^2\geq\left( \dfrac{1}{2}- \dfrac{1}{\theta}\right)m_0=\delta_2 ~~ \forall u\in\mathcal{N}_\Omega.
$$
Now, the result follows by using the equality $I(u)=I(u^{+})+I(u^{-})$ for all $\mathcal{M}_\Omega$.
\fim

\vspace{0.5 cm}

Now, we prove some results related to the following set
$$
\tilde{S}_\lambda:=\{u\in \mathcal{M}_\Omega:\; I(u)<c_\Omega^{*}+\lambda\}.
$$
The above set will be crucial to show the existence of a $(PS)$ sequence of nodal functions associated with $c_\Omega^{*}$.

\begin{lemma}\label{p1.1}
For all $u\in \tilde{S}_\lambda$, we have
$$
0<m_0\leq \|u^{\pm}\|^2\leq\|u\|^2\leq m_\lambda<1,
$$
for $\lambda>0$ sufficiently small.
\end{lemma}
\noindent {\bf Proof.} See Section 4. \fim

\begin{lemma}\label{l1.4}
For each $q>1$, there exists $\delta_{q}>0$ such that
$$
0<\delta_{q}\leq\int_{\Omega}|u^{\pm}|^qdx\leq\int_{\Omega}|u|^qdx, ~~ \forall u\in\tilde{S}_{\lambda}.
$$
\end{lemma}

\noindent {\bf Proof.} ~ See Section 4. \fim

\begin{lemma}\label{l1.6}
There exists $R>0$ such that
$$
I(\dfrac{1}{R}u^{\pm}), \ \ I(Ru^{\pm})<\dfrac{1}{2}I(u^{\pm}), ~\forall u\in\tilde{S}_\lambda
$$

\end{lemma}

\noindent {\bf Proof.} ~ Let $u\in\tilde{S}_\lambda$ and $R>0$. By definition of $I$ and $(f_3)$, 
$$
I\left(\dfrac{1}{R}u^{\pm}\right)=\dfrac{1}{2R^2}\|u^{\pm}\|^2-\int_{\Omega}F\left(\dfrac{1}{R}u^{\pm}\right)dx\leq\dfrac{1}{2R^2}\|u^{\pm}\|^2.
$$
Hence, by Lemma $\ref{p1.1}$
$$
I\left(\dfrac{1}{R}u^{\pm}\right)\leq\dfrac{m_\lambda}{2R^2}.
$$
From this, we can fix $R>0$ large enough such that
$$
\dfrac{m_\lambda}{2R^2}<\delta_2,
$$
which implies, by Corollary $\ref{l1.5}$,
$$
I\left(\dfrac{1}{R}u^{\pm}\right)<\delta_2\leq\dfrac{1}{2}I(u^{\pm}),\ \ \forall u\in\tilde{S}_\lambda.
$$
By $(f_3)$, there are constants $b_1,b_2>0$ verifying
$$
F(t)\geq b_1|t|^{\theta}-b_2,\ \ \forall t\in\mathbb{R},\ \ \forall x\in \Omega.
$$
Then,
$$
I(Ru^{\pm})=\dfrac{R^2}{2}\|u^{\pm}\|^2-\int_{\Omega}F(Ru^{\pm})dx\leq\dfrac{R^2m_\lambda}{2}-b_1R^{\theta}\int_{\Omega}|u^{\pm}|^{\theta}dx+b_2|\Omega|.
$$
By Lemma $\ref{l1.4}$, there is $\delta_\theta>0$ such that
$$
\int_{\Omega}|u^{\pm}|^{\theta}dx\geq\delta_\theta.
$$
Thus,
$$
I(Ru^{\pm})=\dfrac{R^2}{2}\|u^{\pm}\|^2-\int_{\Omega}F(Ru^{\pm})dx\leq\dfrac{R^2m_\lambda}{2}-b_1R^{\theta}\delta_\theta+b_2|\Omega|.
$$
Since $\theta>2$, we conclude that
$$
I(Ru^{\pm})<0<\delta_2\leq\dfrac{1}{2}I(u^{\pm}), \ \ \forall u\in\tilde{S}_\lambda,
$$
for $R>0$ large enough. \fim

\vspace{0.5 cm}

From now on, we consider the following sets

$$
S=\left\{sRu^{+}+tRu^{-}:\ \ u\in \tilde{S}_\lambda\ \ \mbox{and}\ \ s,t\in\left[\dfrac{1}{R^2},1\right]\right\},
$$
$$
P=\{u\in E:\ \ u\geq 0\ \ \mbox{a.e. in}\ \ \Omega\}
$$
and
$$
\Lambda=P\cup(-P).
$$

\begin{lemma}\label{l1.7}
$$
d_0:=dist(S,\Lambda)>0.
$$
\end{lemma}
\noindent {\bf Proof.} The lemma follows by using contradiction
argument combined with Rellich Imbedding. \fim

\begin{proposition}({\bf Main Proposition})\label{p1.2}
Given  $\epsilon,\delta>0$, there exist \linebreak $u\in
I^{-1}([c_\Omega^{*}-2\epsilon,c_\Omega^{*}+2\epsilon])\cap
S_{2\delta}$ verifying
$$
\|I'(u)\|<\dfrac{4\epsilon}{\delta}.
$$

\end{proposition}

\noindent {\bf Proof.}
In fact, otherwise, there exist $\epsilon_o,\delta_o>0$ such that
$$
\|I'(u)\|\geq\dfrac{4\epsilon_o}{\delta_o},\ \ \forall u \in
I^{-1}([c_\Omega^{*}-2\epsilon_o,c_\Omega^{*}+2\epsilon_o])\cap
S_{2\delta_o}.
$$
Thus, for each $n\in\mathbb{N}^*$,
$$\|I'(u)\|\geq\dfrac{4\epsilon_o/n}{\delta_o/n},\ \ \forall u \in
I^{-1}([c_\Omega^{*}-2\epsilon_o,c_\Omega^{*}+2\epsilon_o])\cap
S_{2\delta_o}.
$$
Since
$$
I^{-1}([c_\Omega^{*}-2\epsilon_o/n,c_\Omega^{*}+2\epsilon_o/n])\cap S_{2\delta_o/n}\subset I^{-1}([c_\Omega^{*}-2\epsilon_o,c_\Omega^{*}+2\epsilon_o])\cap S_{2\delta_o},
$$
we get
$$
\|I'(u)\|\geq\dfrac{4\epsilon_o/n}{\delta_o/n},\ \ \forall u \in
I^{-1}([c_\Omega^{*}-2\epsilon_o/n,c_\Omega^{*}+2\epsilon_o/n])\cap
S_{2\delta_o/n}.
$$
Then, we can fix $n\in\mathbb{N}$ large enough such that
\begin{equation}\label{1.10}
\bar{\epsilon}:=\dfrac{\epsilon_o}{n}<\min\left\{\dfrac{2\delta_2}{5},\lambda\right\},\
\ \ \bar{\delta}:=\dfrac{\delta_o}{n}<\dfrac{d_0}{2}
\end{equation}
and
$$
\|I'(u)\|\geq\dfrac{4\bar{\epsilon}}{\bar{\delta}},\ \ \forall u \in
I^{-1}([c_\Omega^{*}-2\bar{\epsilon},c_\Omega^{*}+2\bar{\epsilon}])\cap
S_{2\bar{\delta}}.
$$
The above hypotheses imply that there is continuous map $\eta: E \rightarrow
E$ satisfying:
\begin{enumerate}
\item $\eta(u)=u,\ \ \forall u\notin I^{-1}([c_\Omega^{*}-2\bar{\epsilon},c_\Omega^{*}+2\bar{\epsilon}])\cap S_{2\bar{\delta}}$;

\item $\|\eta(u)-u\| \leq \overline{\delta} ~~ \forall u \in E;$

\item $\eta\left(I^{c_\Omega^{*}+\bar{\epsilon}}\cap S\right)\subset I^{c_\Omega^{*}-\bar{\epsilon}}\cap S_{\bar{\delta}}$;

\item $\eta$ is a homeomorphism.
\end{enumerate}

\vspace{0.5 cm}

From the definition of $c_\Omega^{*}$, for such
$\bar{\epsilon}>0$, there exists ${u}_*\in\mathcal{M}_\Omega$
such that

\begin{equation}\label{1.12}
I({u}_*)<c_\Omega^{*} +\dfrac{\bar{\epsilon}}{2}.
\end{equation}
Now, consider $\gamma:\left[\dfrac{1}{R^2},1\right]\rightarrow
E$ given by
$$
\gamma(s,t)=\eta(sR{u}_*^{+}+tR{u}_*^{-}).
$$
Once ${u}_*^{\pm}\in\mathcal{N}$, 
$$
I(sR{u}_*^{+}+tR{u}_*^{-})=I(sR{u}_*^{+})+I(tR{u}_*^{-})\leq I({u}_*^{+})+({u}_*^{-})=I({u}_*).
$$
Thereby, (\ref{1.10}) and (\ref{1.12}) give 
$$
I(sR{u}_*^{+}+tR{u}_*^{-})\leq I({u}_*)<c_{\Omega}^{*}+\dfrac{\bar{\epsilon}}{2}<c_{\Omega}^{*}+\bar{\epsilon}<c_{\Omega}^{*}+\lambda,
$$
for all $s,t\in\left[\dfrac{1}{R^2},1\right]$. Then, ${u}_* \in\tilde{S}_\lambda$ and
$$
sR{u}_*^{+}+tR{u}_*^{-}\in I^{c_{\Omega}^{*}+\bar{\epsilon}}\cap S,
$$
which implies, by item $3)$,
\begin{equation}\label{1.15}
I(\gamma(s,t))=I(\eta(sR{u_*}^{+}+tR{u_*}^{-}))<c_\Omega^{*}-\bar{\epsilon},\
\ \forall (s,t)\in\left[\dfrac{1}{R^2},1\right]^2.
\end{equation}
From item $2)$,
$$
\|\gamma(s,t)-(sR{u}_*^{+}+tR{u}_*^{-})\|\leq\bar{\delta},
$$
then by the choice of $\bar{\delta}$ made in $(\ref{1.10})$, for $v\in
\Lambda$, we have
$$
\begin{array}{rl}
\|\gamma(s,t)-v\|
&=\|\gamma(s,t)-(sR{u}_*^{+}+tR{u}_*^{-})+(sR{u}_*^{+}+tR{u}_*^{-})-v\|\\[0.5cm]
&\geq\|(sR{u}_*^{+}+tR{u}_*^{-})-v\| - \|\gamma(s,t)-(sR{u}_*^{+}+tR{u}_*^{-})\| \\[0.5cm]
&\geq d_0-\bar{\delta}>d_0-\dfrac{d_0}{2}=d_0>0.
\end{array}
$$
for all $s,t\in\left[\dfrac{1}{R^2},1\right]$. Therefore,
\begin{equation}\label{1.13}
\gamma(s,t)^{\pm}\neq 0,\ \ \forall (s,t)\in
\left[\dfrac{1}{R^2},1\right]^2.
\end{equation}

\begin{claim}\label{af1.1}
There exists $(s_0,t_0)\in \left[\dfrac{1}{R^2},1\right]^2$ such that
$$
I'(\gamma(s_0,t_0)^{\pm})(\gamma(s_0,t_0)^{\pm})=0.
$$
\end{claim}

Suppose, for a moment, that this claim is true. From $(\ref{1.13})$, $\gamma(s_0,t_0)\in\mathcal{M}_\Omega $, and so, 
$$
I(\gamma(s_0,t_0))\geq c_\Omega^{*},
$$
which contradicts $(\ref{1.15})$, proving the proposition. 

\vspace{0.5 cm}

\noindent {\bf Proof of Claim $\ref{af1.1}$:}

\vspace{0.5 cm}

Let us define $ Q:=\left[\dfrac{1}{R^2},1\right]^2 $ and the functions $H, G:Q\rightarrow\mathbb{R}^2$ by
$$
H(s,t):=(I'(\gamma(s,t)^{+}))(\gamma(s,t)^{+}),I'(\gamma(s,t)^{-}))(\gamma(s,t)^{-}))
$$
and
$$
G(s,t):=(I'(sR{u_*}^{+})(sR{u}_*^{+}),I'(tR{u}_*^{-})(tR{u}_*^{-})).
$$
Since
\begin{equation} \label{YYY}
\gamma(s,t)=\eta(sR{u}_*^{+}+tR{u}_*^{-})=sR{u}_*^{+}+tR{u}_*^{-},\ \ \forall (s,t)\in\partial Q,
\end{equation}
we have
$$
\gamma(s,t)^{+}=sR{u}_*^{+}\ \ \mbox{and}\ \ \gamma(s,t)^{-}=tR{u}_*^{-},\ \ \forall (s,t)\in\partial Q,
$$
and $H\equiv G\ \ \mbox{on}\ \ \partial Q.$

To see (\ref{YYY}), let $s=1/R^2$ and $t\in\left[\dfrac{1}{R^2},1\right]$. By Lemma $\ref{l1.6}$,

$$
\begin{array}{rl}
I(sR{u_*}^{+}+tR{u_*}^{-})&=I(\dfrac{1}{R}{u_*}^{+})+I(tR{u_*}^{-})\\[0.5cm]

                                  &<\dfrac{I({u_*}^{+})}{2}+I({u_*}^{-})=I({u_*})-\dfrac{I({u_*}^{+})}{2}.
\end{array}
$$
From $(\ref{1.12})$, Corollary $\ref{l1.5}$ and the choice of
$\bar{\epsilon}>0$ made in $(\ref{1.10})$, we obtain
$$
I(sR{u_*}^{+}+tR{u_*}^{-})<c_\Omega^{*}+\dfrac{\bar{\epsilon}}{2}-\delta_2<c_\Omega^{*}-2\bar{\epsilon},
$$
i.e.,
$$
\dfrac{1}{R}{u_*}^{+}+tR{u_*}^{-}\notin I^{-1}([c_\Omega^{*}-2\bar{\epsilon},c_\Omega^{*}+2\bar{\epsilon}])\cap S_{2\bar{\delta}},
$$
for all $t\in\left[\dfrac{1}{R^2},1\right]$. From this, item 1) yields 
$$
\gamma\left(\dfrac{1}{R^2},t\right)=\eta\left(\dfrac{1}{R}{u_*}^{+}+tR{u_*}^{-}\right)=\dfrac{1}{R}{u_*}^{+}+tR{u_*}^{-}.
$$
The other cases are similar. Then,
$d(H,\dot{Q},(0,0))=d(G,\dot{Q},(0,0))$, but
$d(G,\dot{Q},(0,0))=1\neq 0$. From Brouwer's degree property, there
exists $(s_0,t_0)\in Q$ such that $H(s_0,t_0)=0$, i.e.,
$I'(\gamma(s_0,t_0)^{\pm})(\gamma(s_0,t_0)^{\pm})=0$, and the proof
is complete. \fim


\subsection{Proof of Theorem \ref{bd}}

For each $n\in\mathbb{N}$, consider $\epsilon=\dfrac{1}{4n}$ and $\delta=\dfrac{1}{\sqrt{n}}.$ From Proposition $\ref{p1.2}$, there exists $u_n\in S_{2/\sqrt{n}} $ 
with
$$
u_n\in I^{-1}([c_\Omega^{*}-1/2n,c_\Omega^{*}+1/2n])
$$
and 
$$
\|I'(u_n)\|\leq \dfrac{1}{\sqrt{n}}.
$$
Thus, there is $(v_n)\subset S$ satisfying
$$
I(v_n)\rightarrow c_\Omega^{*}\ \ \mbox{and}\ \ I'(v_n)\rightarrow 0, \ \
$$
in other words, $(v_n)$ is a $(PS)_{c_\Omega^{*}}$ of nodal functions for $I$.

\begin{claim}\label{1.19}
The sequence $(v_n)$ is bounded in $E$ and for a
subsequence of $(v_n)$, still denoted by $(v_n)$,
$$
\limsup_{n \in \mathbb{N}}\|v_n\|^2<1.
$$
\end{claim}

Indeed, since $(v_n)\subset S$, it is easy to see that $(v_n)$ is bounded in $E$. Thus, $I'(v_n)v_n=o_n(1)$ and
$$
c_{\Omega}^{*}+o_n(1)=I(v_n)-\dfrac{1}{\theta}I'(v_n)v_n=\left(\dfrac{1}{2}-\dfrac{1}{\theta}\right)\|v_n\|^2-\int_{\Omega}[F(v_n)-\dfrac{1}{\theta}f(v_n)v_n]dx
$$
The above equality together with $(f_3)$ and Lemma \ref{A} gives
$$
\limsup_n\|v_n\|\leq \frac{c_{\Omega}^{*}}{\Big(\frac{1}{2}-\frac{1}{\theta}\Big)}<1.
$$

Now, let $v_0\in E$ the weak limit of $(v_n)$. Combining 
Claim \ref{1.19} with Lemma \ref{l1.8}, we deduce that $v_0$ is a weak
solution to problem $(P)$. Finally, to conclude the proof, we must prove that
$v_0^{\pm}\neq 0$. We know that

$$
v_n\rightharpoonup v_0\ \ \mbox{in}\ \ H_0^1(\Omega);
$$

$$
v_n(x)\rightarrow v_0(x)\ \ \mbox{a.e. in}\ \ \Omega
$$
and
$$
v_n\rightarrow v_0\ \ \mbox{in}\ \ L^q(\Omega).
$$
On the other hand, using that $v_n\in S$, there are
$s_n,t_n\in\left[\dfrac{1}{R^2},1\right]$ and
$u_n\in\mathcal{M}_\Omega$, such that
$$
v_n=s_nRu_n^{+}+t_nRu_n^{-}\rightharpoonup s_0Ru_0^{+}+t_0Ru_0^{-}\ \ \mbox{in}\ \ E
$$
and
$$
v_n(x)=s_nRu_n^{+}(x)+t_nRu_n^{-}(x)\rightarrow s_0Ru_0^{+}(x)+t_0Ru_0^{-}(x)\ \ \mbox{a.e. in}\ \ \Omega,
$$
for some $s_0, t_0\in\left[\dfrac{1}{R^2},1\right]$, where $u_0\in
E$ is the weak limit of the sequence
$(u_n)\subset\mathcal{M}_{\Omega}$. By uniqueness of limit, we have
$v_0=s_0Ru_0^{+}+t_0Ru_0^{-}$. From Lemma $\ref{l1.4}$, we obtain
$u_0^{\pm}\neq 0$, which implies that $v_0^{+}=s_0Ru_0^{+}\neq 0$
and $v_0^{-}=s_0Ru_0^{-}\neq 0$ and the proof of Theorem \ref{bd} is
complete.


\section{Unbounded Domain}

From now on, we consider the problem $(P)$ with $\Omega= \mathbb{R}^{2}$. From $(V_1)$, it is possible to show that 
$$
\|u\|=\Big(\int_{\R^2}\left(|\nabla u|^2+V(x)|u|^2\right)dx\Big)^{\frac{1}{2}}
$$
is a norm on $H^{1}(\mathbb{R}^{2})$, which is equivalent to the usual norm in $H^{1}(\mathbb{R}^{2})$.  Hereafter,  $E$ denotes $H^{1}(\mathbb{R}^{2})$ endowed with the above norm.

From assumptions $(f_1)$ and $(f_2)$, given $\epsilon>0$, $q\geq 1$
and $\beta>4$, there exists a positive constant
$C=C(\epsilon,q,\beta)$ such that
$$
sf(s),\ \ F(s)\leq \epsilon\dfrac{s^2}{2}+C|s|^{q}\left(e^{\beta\pi s^2}-1\right),\ \ \mbox{for all}\ \ s\in \R.
$$
Thus, by a Trudinger-Moser inequality (\ref{X2}),
we have $F(u)\in L^1(\R^2)$ for all $u\in H^1(\R^2)$. Therefore, the
Euler-Lagrange functional associated with $(P)$ given by
$$
I(u)=\dfrac{1}{2}\|u\|^2-\int_{\R^2}F(u)dx,\ \ u\in E.
$$
is well defined. Furthermore, using standard arguments, we see that
$I$ is a $C^1$ functional on $E$ with
$$
I'(u)v=\int_{\R^2}\left[\nabla u\nabla v+V(x)uv\right]dx-\int_{\R^2}f(u)vdx,\ \ \mbox{for all}\ \ v\in E.
$$
Consequently, critical points of $I$ are precisely the weak solutions
of problem $(P)$. Every nontrivial critical point of $I$ is
contained in the Nehari manifold
$$
\mathcal{N}=\{u\in E\setminus\{0\}:\; I'(u)u=0\}.
$$
A critical point $u\neq 0$ of $I$ is a ground state if $I(u)=c_1$,
where
$$
c_1=\inf_{u\in\mathcal{N}}I(u).
$$
Since we are interested in least energy nodal solution, we define
the nodal Nehari set
$$\mathcal{M}=\{u\in E:\;u^{\pm}\neq 0,\; I'(u^{\pm})u^{\pm}=0\},$$
and
$$c^{*}=\inf_{u\in\mathcal{M}}I(u).$$
Here, it is important to observe that every nodal solution of $(P)$ lies in $\mathcal{M}$.

\vspace{0.5 cm}

Next, we state some necessary results to prove Theorem
\ref{ud}. The proofs of some of them are in Section 4. The first one can be found in Alves, Carri\~ao and Medeiros \cite{ACE}.

\begin{lemma}\label{aux2}
Let $F\in C^2(\R,\R_{+})$ be a convex and even function such that
$F(0)=0$ and $f(s)=F'(s)\geq 0$, $\forall s\in [0,+\infty)$. Then,
for all $t,s\geq 0$
$$
|F(t-s)-F(t)-F(s)|\leq 2(f(t)s+f(s)t).
$$
\end{lemma}

The following two results is essentially due to Alves, do \'O and
Miyagaki and its proof can be found in \cite{AdoOM}.

\begin{theorem}\label{3.1}
Suppose that $(V_1)-(V_2)$ and $(f_1)-(f_5)$ hold. Then
$$
\left\{\begin{array}{lcl} -\Delta u+V_\infty(x)u=f(u),  &\mbox{ in }& \R^2, \\
u\in H^1(\R^2),
\end{array}\right. \leqno (P_\infty)
$$
possesses a positive ground state solution, i. e., there exists
$\bar{u}\in H^1(\R^2)$ such that $\bar{u}>0$,
$I_\infty(\bar{u})=c_\infty$ and $I_\infty'(\bar{u})=0$, where
$$I_\infty(u)=\dfrac{1}{2}\int_{\R^2}\left(|\nabla u|^2+V_\infty(x)u^2\right)dx-\int_{\R^2}F(u)dx,\ \ u\in H^1(\R^2),$$
$$c_\infty=\inf_{u\in\mathcal{N}_\infty}I_\infty(u)$$
and $\mathcal{N}_\infty$ denotes the Nehari manifold
$$\mathcal{N}_\infty=\{u\in H^1(\R^2)\setminus\{0\}:\; I_\infty'(u)u=0\}.$$
\end{theorem}

The second result deal with the asymptotically periodic case.
\begin{theorem}\label{3.2}
Suppose that $(V_1)-(V_2)$ and $(f_1)-(f_5)$ hold. Then,  problem
$(P)$ possesses a positive ground state solution, i. e., there
exists $u_1\in H^1(\R^2)$ such that $u_1>0$, $I(u_1)=c_1$ and
$I'(u_1)=0$.
\end{theorem}

Employing the same arguments explored by  Alves \cite{A}, it is
possible to prove the following result

\begin{theorem}\label{aux1}
Assume that $(f_1)$ and $(f_2)$ hold. Then, any positive solution
$\bar{u}$ of problem $(P_\infty)$ with $\|\bar{u}\|_{H^1(\R^2)}<1$
satifies
$$(I)\ \ \lim_{|x|\rightarrow \infty}\bar{u}(x)=0$$
and
$$(II)\ \ C_1e^{-a|x|}\leq \bar{u}\leq C_2e^{-b|x|}\ \ \mbox{in}\ \ \R^2,$$
where $C_1$ and $C_2$ are positive constants and $0<b<1<a$.
Moreover, we can be chosen $a=1+\delta$, $b=1-\delta$ for
$\delta>0$. The same result hold for $u_1>0$ given in Theorem
\ref{3.2}.
\end{theorem}

The next proposition is a key point in our arguments to get nodal
solution, because it gives an estimate from above of $c^{*}$.
\begin{proposition}\label{cstar}
Suppose that $(V_1)-(V_3)$ hold. Then $c^{*}<c_1+c_\infty$.
\end{proposition}

\noindent {\bf Proof.} See Section 4. \fim

\vspace{0.5 cm}

The below lemma establishes a condition to conclude when the weak limit of a $(PS)$ sequence is nontrivial.

\begin{lemma}\label{neq0}
Assume that $(V_1)-(V_3)$ and $(f_1)-(f_5)$ hold. If $(u_n)\subset
E$ is such that $I(u_n)\rightarrow\sigma$, $u_n\rightharpoonup u$,
$I'(u_n)u_n\rightarrow 0$ and
$$\liminf_{n\rightarrow\infty}\int_{\R^2}f(u_n)f(u_n)dx>0,$$
then $u\neq 0$, provided that $0<\sigma<c_\infty$.
\end{lemma}
\noindent {\bf Proof.} See Section 4. \fim


\subsection{Proof of Theorem \ref{ud}}

Applying Theorem \ref{bd} with
$\Omega=B_n(0)$ and $n \in \mathbb{N}$, there is a nodal 
solution $u_n\in H_0^1(B_n(0))$ for $(P)$ satisfying
$$
I(u_n)=c_n^{*} ~~ \mbox{and} ~~ I'(u_n)=0,
$$
where $c_n^{*}=c_{B_n(0)}^{*}$.
Here, we also denote by $I$ the functional associated with $(P)$, because its restriction to $H_0^1(B_n)$ coincides
with the functional associated with $(P)$.

\begin{claim}\label{cnc} The below limit holds
$$
\lim_{n\rightarrow \infty}c_n^{*}=c^{*}.
$$
\end{claim}

\noindent Indeed, we begin recalling that $(c_n^{*})$ is a non-increasing sequence and
bounded from bellow by $c^{*}$. If $\lim c_n^{*}=\hat{c}>c^{*}$,
then there exists $\phi\in \mathcal{M}$ such that $I(\phi)<\hat{c}$.
Take $(\omega_n)\subset C_0^{\infty}(\R^2)$ and $t_n^{\pm}>0$ such that
$$
\omega_n^{\pm}\neq 0, ~~ \omega_n\rightarrow \phi ~~ \mbox{in} ~~  H^1(\R^2) ~~
\mbox{and} ~~ t_n^{\pm}\omega_n^{\pm}\in \mathcal{N}.
$$
Thereby,
$$
I(\omega_n)=I(\omega_n^{+})+I(\omega_n^{-})\rightarrow I(\phi)\geq c^{*}>0,
$$
$$
I(\omega_n^{\pm})\rightarrow \phi^{\pm},
$$
and
$$
I'(\omega_n^{\pm})\omega_n^{\pm}\rightarrow I'(\phi^{\pm})\phi^{\pm}=0.
$$
Then, if we define
$\phi_n:=t_n^{+}\omega_n^{+}+t_n^{-}\omega_n^{-}\in\mathcal{M}$, by
using similar arguments contained in the proof of Lemma \ref{neq0},
it is possible to prove that
$$
t_n^{\pm}\rightarrow 1 ~~ \mbox{and} ~~ I(t_n^{\pm}\omega_n^{\pm})\rightarrow I(\phi^{\pm}),
$$
leading to,
$$
I(\phi_n)\rightarrow I(\phi).
$$
Therefore, we can fix $n_0\in\mathbb{N}$ such that $I(\phi_{n_0})<\hat{c}, \ \
\forall n\geq n_0$. On the other hand, fixing $n_1 \in \mathbb{N}$ such that $\phi_{n_0}\in \mathcal{M}_{n_1},$ it follows that
$$
c_{n_1}\leq I(\phi_{n_0})<\hat{c},
$$
which contradicts the definition of $\hat{c}$.

From $(f_3)$, we know that $(u_n)$ is a bounded sequence in $E$.
Thus, we can assume that $(u_n)$ is weakly convergent to $u$, for
some $u\in E$. Once 
$$
c^{*}=\lim_nc_n^{*}=\lim_nI(u_n)
$$
and
$$
I'(u_n)v=0,\ \ \mbox{for all}\ \ v\in H_0^1(B_n),
$$
a direct computation gives that $u$ is a weak solution for $(P)$. Now, our goal is proving that
$$
u\in\mathcal{M} ~~ \mbox{and} ~~ I(u)=c^{*}.
$$
In fact, taking a subsequence if necessary, we can assume that
$$
I(u_n^{\pm})\rightarrow\sigma^{\pm},\ \ \mbox{where}\ \ c^{*}=\sigma^{+}+\sigma^{-}.
$$

Using that $u_n^{+},u_n^{-}\in\mathcal{N}$, we derive $\sigma^{\pm}\geq
c_1>0$.
 From Proposition \ref{cstar}, it follows that $\sigma^{\pm}<c_\infty$.
Since
$$
\liminf_{n\rightarrow\infty}\int_{\R^2}f(u_n^{\pm})u_n^{\pm}>0,
$$
Lemma \ref{neq0} yields $u^{\pm}\neq 0$. Therefore, $u\in
\mathcal{M}$ and $I(u)\geq c^{*}$. To complete the proof, by Fatou's
Lemma, we see that
$$2c^{*}=\liminf_{n\rightarrow\infty}\left[2I(u_n)-I'(u_n)u_n\right]=\liminf_{n\rightarrow \infty}\int_{\R^2}\left(f(u_n)u_n-2F(u_n)\right)dx$$
$$\geq\int_{\R^2}\left(f(u)u-2F(u)\right)dx=2I(u)-I'(u)u=2I(u)\geq 2c^{*}.$$
Hence, $I(u)=c^{*}$, which proves that $(P)$ has a nodal
solution. In order to establish that the nodal solution has
exactly two nodal domain, we refer the reader to \cite[Theorem 2.3]{BW}.  \fim

\section{Nonexistence result}

In this section, we prove a nonexistence result of least energy
nodal solution for the following autonomous problem
$$
\left\{
\begin{array}{l}
-\Delta u+V_0u=f(u),  \mbox{ in } \R^2, \\[0.3cm]
u\in H^1(\R^2),
\end{array}
\right. \eqno (Q)
$$
that is, we prove that
$$
\hat{c}:=\inf_{\mathcal{M}} J
$$
is not attained, where $J$ is the energy functional defined on $H^1(\R^2)$ associated with $(Q)$ and $\mathcal{M}$ is the nodal Nehari set
$$
\mathcal{M}:=\{u\in H^1(\R^2):~~ u^{\pm}\neq 0~~\mbox{and}~~J'(u^{\pm})u^{\pm}=0\}.
$$
 For this, we define
$$
f_{+}(t)=
\left\{
\begin{array}{l}
 f(t), ~ ~ t\geq 0, \\
 0, ~~ t \leq 0
\end{array}
\right.
$$
and the functional $J_{+}$ defined on $H^1(\R^2)$ by
$$
J_{+}(u):=\int_{\R^2}(|\nabla u|^2+V_0|u|^2)dx-\int_{\R^2}F_{+}(u)dx,
$$
where $F_{+}$ is the primitive of $f_{+}$ with $F_{+}(0)=0$. From  \cite[Theorem 1.1]{AdoOM}, the below number
$$
c_{+}=\inf_{\mathcal{N}_{+}} J_{+}
$$
where
$$
\mathcal{N}_{+}:=\{u\in H^1(\R^2)\setminus\{0\}:~~J_{+}'(u)u=0\},
$$
is a critical value of $J_{+}$. Let $v$ be the corresponding critical point. It is easy to see that $v^{-}=0$. Thus, $v$ is nonnegative and by the maximum principle, $v>0$ on $\R^2$. In particular, $v$ is a positive critical point of $J$.

Analogously, if we define
$$
f_{-}(t)=\left\{
\begin{array}{l}
0, ~ ~ t\geq 0, \\
f(t), ~ ~ ~ t<0,
\end{array}
\right.
$$
and denote by $J_{-}$ the corresponding functional and by $\mathcal{N}_{-}$ the Nehari manifold, then
$$
c_{-}:=\inf_{\mathcal{N}_{-}} J_{-}
$$
is a critical value of $J_{-}$.

The next proposition is a key point in our argument to prove the nonexistence result, because it gives an exact estimate of $\hat{c}$.

\begin{proposition}\label{cc}
Under assumptions $(f_1)-(f_5)$, we have
$$
\hat{c}=c_{+}+c_{-}.
$$
\end{proposition}

\noindent {\bf Proof.} ~ Let $v,w\in H^1(\R^2)$ verifying

$$
J_{+}(v)=c_{+},~~  J'_{+}(v)=0,~~ v(x)>0, ~~\forall x\in\R^2,
$$

$$
J_{-}(w)=c_{-},~~  J'_{-}(w)=0~~ w(x)<0, ~~\forall x\in\R^2
$$
and consider the functions
$$
v_R(x):=\varphi\left(\frac{x}{R}\right)v(x)~~~\mbox{and}~~~ w_{R,n}:=\varphi\left(\frac{x-x_n}{R}\right)w(x-x_n),
$$
where $\varphi\in C^{\infty}_0(\R^2)$ is a cut-off function satisfying
$$
supp~\varphi\subset B_2(0), ~0\leq\varphi\leq 1, ~\varphi=1 ~ \mbox{on} ~~ B_1(0) ~~ \mbox{and} ~~ x_n=(n,0).
$$
Clearly, for $n$ large enough,
$$
supp ~ v_R\cap supp ~ w_{R,n}=\emptyset.
$$
Let $t_R,s_R$ be the positive real numbers such that
$$
J'(t_Rv_R)t_Rv_R=0~~~ \mbox{and}~~~ J'(s_Rw_{R,n})s_Rw_{R,n}=0.
$$
Since
$$
t_R^2\int_{\R^2}\left(|\nabla v_R|^2+V_0|v_R|^2\right)dx=\int_{\R^2}f_{+}(t_Rv_R)t_Rv_R
$$
and $v_R\rightarrow v$ in $H^1(\R^2)$ as $R\rightarrow +\infty$, it is possible to show, by using similar arguments given in the proof of Lemma \ref{neq0}, that
$t_R\rightarrow 1$, as $R\rightarrow +\infty$. Similarly,
$$
s_R^2\int_{\R^2}\left(|\nabla w_{R,n}|^2+V_0|w_{R,n}|^2\right)dx=\int_{\R^2}f_{+}(s_Rw_{R,n})s_Rw_{R,n}.
$$
Since $w_R\rightarrow w$ in $H^1(\R^2)$ as $R\rightarrow +\infty$, we derive that
$s_R\rightarrow 1$, as $R\rightarrow +\infty$.
Now, note that $u_R:=t_Rv_R+s_Rw_{R,n}\in\mathcal{M}$ with
$$
u_R^{+}=t_Rv_R~~~\mbox{and}~~~u_R^{-}=s_Rw_{R,n}
$$
for $n\in\mathbb{N}$ large enough. Then,
$$
\hat{c}\leq J(t_Rv_R+s_Rw_{R,n})=J(t_Rv_R)+J(s_Rw_{R,n})
$$
Using the invariance of $\R^2$ under translations, we obtain by taking $R \to +\infty$
$$
\hat{c}\leq J(v)+J(w). 
$$
Since $J(v)=J_{+}(v)=c_{+}$ and $J(w)=J_{-}(w)=c_{-}$, it follows that
$$
\hat{c}\leq c_{+}+c_{-}.
$$
On the other hand, it is obvious that $\hat{c}\geq c_{+}+c_{-}$. Therefore, we can conclude that $\hat{c}= c_{+}+c_{-}$.
\fim

\vspace{0.5 cm}

\noindent {\bf Proof of Theorem \ref{ud2}.} Suppose by contradiction that there exists $u\in\mathcal{M}$ such that $J(u)=\hat{c}$. Thus, $u^{+}\in\mathcal{N}_{+}$ and $u^{-}\in\mathcal{N}_{-}$, from where it follows that
$$
c_{+}+c_{-}\leq J_{+}(u^{+})+J_{-}(u^{-})=J(u)=\hat{c}=c_{+}+c_{-},
$$
and so, 
$$
J_{+}(u^{+})=c_{+}~~\mbox{and}~~ J_{-}(u^{-})=c_{-}.
$$
Thereby, $u^{+}$ is a critical point of $J_{+}$ and $u^{-}$ is a critical point of $J_{-}$.Then, by maximum principle, we must have 
$$
u^{+}(x)>0,~~\mbox{for all}~~x\in\R^2
$$
and
$$
u^{-}(x)<0,~~\mbox{for all}~~x\in\R^2,
$$
which is impossible.
\fim

\begin{remark}
A version of Theorem \ref{ud2} can be make for $N \geq 3$, by supposing that $f$ has a subcritical growth.
\end{remark}

\begin{remark}
We can define $H_{r}^{1}(\R^2):=\{u\in H: ~ u ~ \mbox{is a radial function}\}$, $\mathcal{M}_r:=\mathcal{M}\cap H_{r}^{1}(\R^2)$ and $c_r^{*}=\inf_{\mathcal{M}_r} J$.
Under the assumptions of Theorem \ref{ud2}, there exist a minimizer $u\in\mathcal{M}_r$ which is a critical point of $I$ on $H^1(\R^2)$. To prove this, we combine the symmetric criticality principle with arguments as in the proof of Theorem \ref{bd}. It is clear that $c^{*}\leq c_r^{*}$, and so, as a consequence of our nonexistence result, we have $c^{*}< c_r^{*}$. A similar inequality in bounded domain like annulus for $N \geq 3$ was proved in \cite{BWW} .
\end{remark}

\section{Proof of lemmas and propositions}

\noindent {\bf Proof of Lemma \ref{l1.8}.} From $(f_1)$,
$$
|f(u_n)u_n|\leq C|u_n|e^{4 \pi|u_n|^2}, ~~\forall n\in\mathbb{N}.
$$
We claim that
$$
\int_{\Omega}|u_n|e^{4 \pi |u_n|^2}dx\rightarrow
\int_{\Omega}|u|e^{4 \pi |u|^2}dx,~~\mbox{as}~~n\rightarrow\infty.
$$
Effectively, consider $t>1$ with $t\approx 1$. Note that
$$
\int_{\Omega}\left(e^{4 \pi|u_n|^2}\right)^tdx=\int_{\Omega}e^{4 \pi
t\|u_n\|^2\left(\frac{|u_n|}{\|u_n\|}\right)^2}dx\leq\int_{\Omega}e^{4 \pi
tb\left(\frac{|u_n|}{\|u_n\|}\right)^2}dx.
$$

Now, since $b<1$, we can fix $t>1$ with $t\approx 1$, such that $t b<1$.
Consequently, by Trudinger-Moser inequality,
$$
\sup_n\int_{\Omega}\left(e^{4 \pi|u_n|^2}\right)^tdx\leq\sup_{\|v\|\leq
1}\int_{\Omega}e^{4 \pi tb|v|^2}dx<\infty.
$$
Thus, the sequence $(e^{4 \pi |u_n|^2})$ is bounded in $L^t(\Omega)$
and
$$
e^{4 \pi |u_n(x)|^2}\rightarrow e^{4 \pi |u(x)|^2}~~\mbox{a.e. in}~~
\Omega.
$$
This implies that,
\begin{equation}\label{cb}
e^{4 \pi |u_n|^2}\rightharpoonup e^{4 \pi |u|^2}~~\mbox{in}~~
L^t(\Omega).
\end{equation}
On the other hand,
\begin{equation}\label{ca}
|u_n| \rightarrow |u| ~~\mbox{in} ~~ L^{t'}(\Omega)
\end{equation}
where $1/t+1/t'=1$. Now, (\ref{cb}) combined with (\ref{ca}) gives
$$
\int_{\Omega}|u_n|e^{4 \pi |u_n|^2}dx\rightarrow
\int_{\Omega}|u|e^{4 \pi |u|^2}dx.
$$
Hence,
$$
|u_n|e^{4 \pi |u_n|^2} \to |u|e^{4 \pi |u|^2} ~~ \mbox{in} ~~ L^{1}(\Omega).
$$
Thus, for some subsequence, there is $h \in L^{1}(\Omega)$ verifying
$$
|u_n|e^{4 \pi |u_n|^2} \leq h ~~ \mbox{a.e. in } ~~ \Omega.
$$
Thereby,
$$
|f(u_n)u_n| \leq h ~~ \mbox{a.e. in } ~~ \Omega.
$$
By Lebesgue's Theorem, it follows that
$$
\int_{\Omega}f(u_n)u_ndx\to \int_{\Omega}f(u)udx.
$$
The proof of (\ref{1.23}) follows by using the same type of arguments. \fim

\vspace{0.5 cm}

\noindent {\bf Proof of Lemma \ref{p1.1}.} Since
$\tilde{S}_\lambda\subset\mathcal{M}$, in view of Corollary
\ref{c1.1}, we only need to prove that there exist $m_\lambda>0$
such that
$$
\|u\|^2\leq m_\lambda<1 ~~ \forall u \in \tilde{S}_{\lambda}.
$$
For each $u\in \tilde{S}_{\lambda}$, we have
$$
c_{\Omega}^{*} +\lambda\geq I(u)=I(u)-\dfrac{1}{\theta}I'(u)u=\left(\dfrac{1}{2}-\dfrac{1}{\theta}\right)\|u\|^2-\int_{\Omega}\left(F(x,u)-\dfrac{1}{\theta}f(x,u)u\right)dx.
$$
From Ambrosetti-Rabinowitz condition $(f_3)$, 
$$
c_{\Omega}^{*} +\lambda\geq\left(\dfrac{1}{2}-\dfrac{1}{\theta}\right)\|u\|^2.
$$
On the other hand, by Lemma \ref{A}, we can fix $\lambda>0$
sufficiently small such that
$$
A +\lambda<\left(\dfrac{1}{2}-\dfrac{1}{\theta}\right),
$$
where $A$ was given in (\ref{A1}). Therefore
$$
\|u\|^2\leq \dfrac{c_{\Omega}^{*} +\lambda}{\left(\dfrac{1}{2}-\dfrac{1}{\theta}\right)}\leq m_{\lambda}<1,
$$
where
$$
m_\lambda:=\dfrac{A+\lambda}{\left(\dfrac{1}{2}-\dfrac{1}{\theta}\right)}.
$$
\fim

\vspace{0.5 cm}

\noindent{\bf Proof of Lemma \ref{l1.4}.} Since $u\in
\tilde{S}_\lambda\subset\mathcal{M}$,
$$
\|u^{\pm}\|^2=\int_{\Omega}f(u^{\pm})u^{\pm}dx.
$$
Then, from $(f_1)$,
$$
\|u^{\pm}\|^2\leq C\int_{\Omega}|u^{\pm}|e^{4 \pi|u^{\pm}|^2}dx.
$$
Using Sobolev imbedding and H\"{o}der inequality, for $1<t_1$ and
$1<t_2\approx 1$ such that $1/t_1+1/t_2=1$, we obtain
$$
\|u^{\pm}\|^2\leq |u^{\pm}|_{L^{t_1}}\left(\int_{\Omega}e^{4 \pi t_2|u^{\pm}|^2}dx\right)^{1/t_2}.
$$
From Corollary $\ref{c1.1}$,
$$
m_0\leq |u^{\pm}|_{L^{t_1}}\left(\int_{\Omega}e^{4 \pi t_2\|u^{\pm}\|^2\left(\frac{|u^{\pm}|}{\|u^{\pm}\|}\right)^2}dx\right)^{1/t_2},
$$
and by Lemma \ref{p1.1}, it follows that
$$
m_0\leq |u^{\pm}|_{L^{t_1}}\left(\int_{\Omega}e^{4 \pi t_2m_\lambda\left(\frac{|u^{\pm}|}{\|u^{\pm}\|}\right)^2}dx\right)^{1/t_2}.
$$
Since $m_\lambda< 1$, we can fix $1<t_2$ near $1$ such
that $ t_2 m_\lambda<1$. From Trudinger-Moser inequality (\ref{X1}),
there exists a constant $C>0$ such that
$$
\int_{\Omega}e^{4 \pi t_2m_\lambda\left(\frac{|u^{\pm}|}{\|u^{\pm}\|}\right)^2}dx\leq C \,\,\, \forall u\in\tilde{S}_\lambda.
$$
Thereby, for some $C_1>0$,
$$
C_1 \leq |u^{\pm}|_{L^{t_1}} ~~ \forall ~~ u\in\tilde{S}_\lambda.
$$
Now, the lemma follows applying interpolation. \fim

\vspace{0.5 cm}

\noindent {\bf Proof of Proposition \ref{cstar}.} Let $\bar{u}$ be a ground
state solution of $(P_\infty)$ and $u_1$ is a positive ground state
of $(P)$ given by Theorems  \ref{3.1} and \ref{3.2},
respectively. Let us define $\bar{u}_n(x)=\bar{u}(x-x_n)$, where
$x_n=(0,n)$ and for $\alpha,\beta>0$
$$
\begin{array}{l}
h^{\pm}(\alpha,\beta,n)=\displaystyle \int_{\R^2}\left(|\nabla(\alpha u_1-\beta\bar{u}_n)^{\pm}|^2+V(x)|(\alpha u_1-\beta\bar{u}_n)^{\pm}|^2\right)dx\\
\mbox{}\\
-\displaystyle  \int_{\R^2}f((\alpha u_1-\beta\bar{u}_n)^{\pm})(\alpha u_1-\beta\bar{u}_n)^{\pm}dx.
\end{array}
$$
Recalling that $I'(u_1)u_1=0$ and using $(f_4)$, we get
$$
\int_{\R^2}\left(|\nabla
(u_1/2)|^2+V(x)(u_1/2)^2\right)dx-\int_{\R^2}f(u_1/2)(u_1/2)
$$
\begin{equation}\label{4.3}
=\int_{\R^2}\left(\dfrac{f(u_1)}{u_1}-\dfrac{f(u_1/2)}{(u_1/2)}\right)\left(\dfrac{u_1}{2}\right)^2dx>0.
\end{equation}
and
$$
\int_{\R^2}\left(|2\nabla
(u_1)|^2+V(x)|2u_1|^2\right)dx-\int_{\R^2}f(2u_1)(2u_1)
$$
\begin{equation}\label{4.4}
=\int_{\R^2}\left(\dfrac{f(u_1)}{u_1}-\dfrac{f(2u_1)}{2u_1}\right)\left(2u_1\right)^2dx<0.
\end{equation}
By $(V_2)$, for $n$ large enough there holds
\begin{equation}\label{4.5}
\int_{\R^2}\left(|\nabla
(\bar{u}_n/2)|^2+V(x)(\bar{u}_n/2)^2\right)dx-\int_{\R^2}f(\bar{u}_n/2)(\bar{u}_n/2)>0
\end{equation}
and
\begin{equation}\label{4.6}
\int_{\R^2}\left(|\nabla
(2\bar{u}_n)|^2+V(x)(2\bar{u}_n)^2\right)dx-\int_{\R^2}f(2\bar{u}_n)(2\bar{u}_n)<0.
\end{equation}
Hence, from (\ref{4.3})-(\ref{4.6}),  there exists $n_0>0$ such
that

\begin{equation}\label{5.4}
\left\{
\begin{array}{l}
h^{+}(1/2,\beta,n)>0,\\
\mbox{}\\
h^{+}(2,\beta,n)<0
\end{array}
\right.
\end{equation}
for $n\geq n_0$ and $\beta\in[1/2,2]$. Now, for all
$\alpha\in[1/2,2]$ we have

\begin{equation}\label{5.5}
\left\{
\begin{array}{l}
h^{-}(\alpha,1/2,n)>0,\\
h^{-}(\alpha,2,n)<0.
\end{array}
\right.
\end{equation}
From this, we can apply a variant of the Mean Value Theorem due to
Miranda \cite{M}, to obtain $\alpha^{*},\beta^{*}\in[1/2,2]$ such
that $h^{\pm}(\alpha^{*},\beta^{*},n)=0$, for any $n\geq n_0$. Thus,
$$
\alpha^{*}u_1+\beta^{*}\bar{u}_n\in\mathcal{M},\ \ \mbox{for}\ \ n\geq n_0.
$$

In view of the definition of $c^{*}$, it suffices to show that
$$
\sup_{\frac{1}{2}\leq\alpha,\beta\leq 2}I(\alpha u_1+\beta\bar{u}_n)<c_1+c_\infty\ \ \mbox{for}\ \ n\geq n_0.
$$
In order to do this, first we use Lemma \ref{aux2} to get the ensuing estimate
$$
I(\alpha u_1-\beta\bar{u}_n)\leq\dfrac{1}{2}\int_{\R^2}(|\nabla (\alpha u_1))|^2+|\nabla (\beta\bar{u}_n))|^2)dx+\dfrac{1}{2}\int_{\R^2}V(x)(
|\alpha u_1|^2+|\beta\bar{u}_n|^2)dx
$$
$$
-\alpha\beta\int_{\R^2}(\nabla u_1\nabla\bar{u}_n+V(x)u_1\bar{u}_n)dx-A_1,
$$
where
$$
A_1=\int_{\R^2}F(\alpha u_1)dx+\int_{\R^2}F(\beta\bar{u}_n)dx-2\int_{\R^2}\left[f(\alpha u_1)\beta\bar{u}_n+f(\beta\bar{u}_n)\alpha u_1\right]dx
$$
Since $u_1$ is a positive solution of $(P)$, we know that
$$
\int_{\R^2}(\nabla u_1\nabla\bar{u}_n+V(x)u_1\bar{u}_n)dx\geq 0.
$$
Therefore
\begin{equation}\label{4.7}
I(\alpha u_1-\beta\bar{u}_n)\leq I(\alpha
u_1)+I_\infty(\beta\bar{u}_n)+2\alpha\int_{\R^2}f(\beta\bar{u}_n)
u_1dx+2\beta\int_{\R^2}f(\alpha u_1)\bar{u}_ndx
\end{equation}
$$
+\dfrac{\beta^2}{2}\int_{\R^2}(V(x)-V_\infty(x))\bar{u}_n^2dx.
$$
From $(V_3)$,
$$
\int_{\R^2}(V(x)-V_\infty(x))\bar{u}_n^2dx\leq-Ce^{-\mu n}
$$
and by $(f_1)-(f_2)$,
$$
\int_{\R^2}f(\alpha u_1)\bar{u}_ndx\leq\epsilon\alpha\int_{\R^2}u_1\bar{u}_ndx+C\int_{\R^2}\left(e^{\tau\alpha^2u_1^2}-1\right)u_1 \bar{u}_ndx,\ \ \mbox{for}\ \ \tau>4\pi.
$$
Notice that from Theorem \ref{aux1}, 
$$
\int_{B_{n/2}}u_1\bar{u}_ndx\leq C_2\int_{B_{n/2}(0)}u_1e^{-b|x-x_n|}dx.
$$
Once $|x-x_n|\geq |x_n|-|x|=n-|x|$ and $|x|\leq n/2$, we find that 
$|x-x_n|\geq n/2$, from where it follows that 
$$
\int_{B_{n/2}}u_1\bar{u}_ndx\leq C_2\int_{B_{n/2}}u_1e^{-bn/2}dx=Ce^{-bn/2}
$$
and
$$
\int_{\R^2\setminus B_{n/2}}u_1\bar{u}_ndx\leq C_2\int_{\R^2\setminus B_{n/2}}e^{-b|x|}\bar{u}_ndx\leq C_2e^{-bn/2}\int_{\R^2}\bar{u}_ndx=C_2e^{-bn/2}\int_{\R^2}\bar{u}dx.
$$
Therefore
$$
\int_{\R^2}u_1\bar{u}_ndx\leq Ce^{-bn/2}.
$$
Moreover, since $u_1 \in L^{\infty}(\mathbb{R}^{2})$,
$$
\int_{\R^2}\left(e^{\tau\alpha^2u_1^2}-1\right)u_1\bar{u}_ndx\leq C\int_{\R^2}u_1\bar{u}_n\leq Ce^{-bn/2}.
$$
Therefore
$$
\int_{\R^2}f(\alpha u_1)\bar{u}_ndx\leq Ce^{-bn/2} ~~ \mbox{and} ~~ \int_{\R^2}f(\beta\bar{u}_n)u_1dx\leq Ce^{-bn/2}.
$$
Then, from (\ref{4.7}),
$$
I(\alpha u_1-\beta\bar{u}_n)\leq\sup_{\alpha\geq 0}I(\alpha u_1)+\sup_{\beta\geq 0}I(\beta\bar{u}_n)+C(e^{-bn/2}-e^{-\mu n}).
$$
Since $\mu<1/2$, for $n$ large enough, we know that
$$
e^{-bn/2}-e^{-\mu n} <0,
$$
from where it follows that
$$
\sup_{1/2\leq\alpha,\beta\leq 2}I(\alpha u_1-\beta \bar{u}_n)<c_1+c_\infty.
$$
Consequently
$$
c^{*}<c_1+c_\infty,
$$
finishing the proof of the proposition. \fim


\vspace{0.5 cm}

\noindent {\bf Proof of Lemma \ref{neq0}.} Suppose, by contradiction that
$u\equiv 0$. From $(V_2)$, given $\epsilon>0$ there exists
$R=R(\epsilon)>0$ such that
$$
|V(x)-V_\infty(x)|<\epsilon,\ \ \mbox{for} ~~~~ |x|\geq R.
$$
As a consequence of $u\equiv 0$, we get
$$
\int_{B_R}|V(x)-V_\infty(x)||u_n|^2dx \to 0.
$$
The below inequality 
$$
\int_{\R^2}|V(x)-V_\infty(x)||u_n|^2dx\leq\int_{B_R}|V(x)-V_\infty(x)||u_n|^2dx+\epsilon\int_{\R^2\setminus B_R}|u_n|^2dx,
$$
together with the boundedness of $(u_n)$ in $H^{1}(\mathbb{R}^{2})$ yields
$$
|I(u_n)-I_\infty(u_n)|\rightarrow 0\ \ \mbox{as}\ \ n\rightarrow \infty.
$$
A similar argument shows that
$$|I'(u_n)u_n-I_\infty'(u_n)u_n|\rightarrow 0\ \ \mbox{as}\ \ n\rightarrow \infty.$$
Consequently,
\begin{equation}\label{3.18}
I_\infty(u_n)=\sigma+o_n(1)\ \ \mbox{and}\ \
I_\infty'(u_n)u_n=o_n(1).
\end{equation}
In what follows, we fix $s_n>0$ verifying
$$
s_nu_n\in\mathcal{N}_\infty.
$$
We claim that $(s_n)$ converges to $1$ as $n\rightarrow\infty$.
Effectivelly, we start proving that
\begin{equation}\label{3.20}
\limsup s_n\leq 1.
\end{equation}
Suppose by contradiction that there exists a subsequence of $(s_n)$,
still denoted by $(s_n)$, such that $s_n\geq 1+\delta$ for all
$n\in\mathcal{N}$, for some $\delta>0$. From (\ref{3.18}),
\begin{equation}\label{3.20.1}
\int_{\R^2}\left(|\nabla
u_n|^2+V_\infty(x)|u_n|^2\right)dx=\int_{\R^2}f(u_n)u_ndx+o_n(1)
\end{equation}
On the other hand, since $s_nu_n\in\mathcal{N}_\infty$,
$$s_n\int_{\R^2}\left(|\nabla
u_n|^2+V_\infty(x)|u_n|^2\right)dx=\int_{\R^2}f(s_nu_n)u_ndx.$$
Consequently
\begin{equation}\label{3.21}
\int_{\R^2}\left(\dfrac{f(s_nu_n)}{s_nu_n}-\dfrac{f(u_n)}{u_n}\right)|u_n|^2dx=o_n(1).
\end{equation}

We claim that there exist $(y_n)\subset\mathbb{Z}^{2}$ with
$|y_n|\rightarrow\infty$, $r>0$ and $\beta>0$ such that
$$\int_{B_r(y_n)}u_n^2dx\geq\beta>0.$$
Indeed, contrary case, using a version of
Lions' results to critical growth in $\R^2$ due to Alves, do \'O and
Miyagaki \cite{AdoOM}, we derive
$$
\lim_{n \to +\infty}\int_{\R^2}f(u_n)u_ndx=0,
$$
which is contrary to our assumption.

Now, let $v_n(x):=u_n(x+y_n)$. Once that $(u_n)$ is bounded in $H^{1}(\mathbb{R}^{2})$, it is easy to show that $(v_n)$ is also bounded in $H^{1}(\mathbb{R}^{2})$. Therefore, for some subsequence, we can assume that $(v_n)$ is weakly convergent, and we will denote by $\tilde{v}$ its weak limit in $H^{1}(\mathbb{R}^{2})$. Observing that
$$
\int_{B_r(0)}|v_n|^2dx=\int_{B_r(y_n)}|u_n|^2dx\geq\beta>0,
$$
we deduce that $\tilde{v}\neq 0$ in $H^{1}(\mathbb{R}^{2})$. Now, (\ref{3.21}), $(f_4)$ and Fatou's Lemma load to
$$
0<\int_{\R^2}\left(\dfrac{f((1+\delta)\tilde{v})}{(1+\delta)\tilde{v}}-\dfrac{f(\tilde{v})}{\tilde{v}}\right)\tilde{v}^2dx\leq 0,
$$
which is impossible. Hence
$$
\limsup_{n\rightarrow \infty} s_n\leq 1.
$$
If $s_0=\displaystyle \limsup_{n\rightarrow \infty} s_n<1$, we can assume that $s_n<1$ for $n$ large enough. Then, by Fatou's
Lemma
$$0<\int_{\R^2}\left(\dfrac{f(\tilde{v})}{\tilde{v}}-\dfrac{f(s_o\tilde{v})}{s_o\tilde{v}}\right)\tilde{v}^2dx\leq 0\ \ \mbox{if}\ \  s_o>0$$
and
$$0<\int_{\R^2}f(\tilde{v})\tilde{v}\leq 0\ \ \mbox{if}\ \ s_o=0,$$
which are impossible. Hence, $\displaystyle \limsup_{n \to \infty} s_n=1$, and so, for some subsequence,
\begin{equation}\label{3.22}
\lim_{n \to \infty} s_n=1.
\end{equation}

As a consequence of (\ref{3.22}),
$$\int_{\R^2}F(s_nu_n)dx-\int_{\R^2}F(u_n)dx=o_n(1)$$
and
$$
(s_n^2-1)\int_{\R^2}\left(|\nabla u_n|^2+V_\infty(x)|u_n|^2\right)dx=o_n(1),
$$
from where it follows that
$$
I_\infty(s_nu_n)=I_\infty(u_n)+o_n(1).
$$
Then
$$
c_\infty\leq I_\infty(s_nu_n)=\sigma+o_n(1).
$$
Taking $n\rightarrow+\infty$, we find $c_\infty\leq\sigma$, which is
impossible because $\sigma<c_\infty$. This contradiction comes from
the assumption that $u\equiv 0.$ \fim

\end{document}